\documentclass[11pt, a4paper]{article}
\usepackage[cp1251]{inputenc}
\usepackage[russian]{babel}
\usepackage{amsmath} \usepackage{euscript}
\usepackage{mymatrix2}
\usepackage{mymatrix}
\usepackage{longtable}
\usepackage{amssymb}
\begin{document}

\newcounter{bnomer} \newcounter{snomer}
\newcounter{bsnomer}
\setcounter{bnomer}{0}
\renewcommand{\thesnomer}{\thebnomer.\arabic{snomer}}
\renewcommand{\thebsnomer}{\thebnomer.\arabic{bsnomer}}
\renewcommand{\refname}{\begin{center}\large{\textbf{References}}\end{center}}

\newcommand{\sect}[1]{%
\setcounter{snomer}{0}\setcounter{bsnomer}{0}
\refstepcounter{bnomer}
\par\bigskip\begin{center}\large{\textbf{\arabic{bnomer}. {#1}}}\end{center}}
\newcommand{\sst}{%
\refstepcounter{bsnomer}
\par\bigskip\textbf{\arabic{bnomer}.\arabic{bsnomer}. }}
\newcommand{\defi}[1]{%
\refstepcounter{snomer}
\par\medskip\textbf{Definition \arabic{bnomer}.\arabic{snomer}. }{#1}\par\medskip}
\newcommand{\theo}[2]{%
\refstepcounter{snomer}
\par\textbf{Теорема \arabic{bnomer}.\arabic{snomer}. }{#2} {\emph{#1}}\hspace{\fill}$\square$\par}
\newcommand{\mtheo}[1]{%
\refstepcounter{snomer}
\par\medskip\textbf{Theorem \arabic{bnomer}.\arabic{snomer}. }{\emph{#1}}\par\medskip}
\newcommand{\theobp}[2]{%
\refstepcounter{snomer}
\par\textbf{Теорема \arabic{bnomer}.\arabic{snomer}. }{#2} {\emph{#1}}\par}
\newcommand{\theop}[2]{%
\refstepcounter{snomer}
\par\textbf{Теорема \arabic{bnomer}.\arabic{snomer}. }{\emph{#1}}
\par\textsc{Proof}. {#2}\hspace{\fill}$\square$\par}
\newcommand{\theosp}[2]{%
\refstepcounter{snomer}
\par\textbf{Теорема \arabic{bnomer}.\arabic{snomer}. }{\emph{#1}}
\par\textbf{Схема доказательства}. {#2}\hspace{\fill}$\square$\par}
\newcommand{\exam}[1]{%
\refstepcounter{snomer}
\par\medskip\textbf{Example \arabic{bnomer}.\arabic{snomer}. }{#1}\par\medskip}
\newcommand{\deno}[1]{%
\refstepcounter{snomer}
\par\textbf{Definition \arabic{bnomer}.\arabic{snomer}. }{#1}\par}
\newcommand{\post}[1]{%
\refstepcounter{snomer}
\par\textbf{Предложение \arabic{bnomer}.\arabic{snomer}. }{\emph{#1}}\hspace{\fill}$\square$\par}
\newcommand{\postp}[2]{%
\refstepcounter{snomer}
\par\medskip\textbf{Proposition \arabic{bnomer}.\arabic{snomer}. }{\emph{#1}}%
\ifhmode\par\fi\textsc{Proof}. {#2}\hspace{\fill}$\square$\par\medskip}
\newcommand{\lemm}[1]{%
\refstepcounter{snomer}
\par\textbf{Lemma \arabic{bnomer}.\arabic{snomer}. }{\emph{#1}}\hspace{\fill}$\square$\par}
\newcommand{\lemmp}[2]{%
\refstepcounter{snomer}
\par\medskip\textbf{Lemma \arabic{bnomer}.\arabic{snomer}. }{\emph{#1}}
\par\textsc{Proof}. {#2}\hspace{\fill}$\square$\par\medskip}
\newcommand{\coro}[1]{%
\refstepcounter{snomer}
\par\textbf{Следствие \arabic{bnomer}.\arabic{snomer}. }{\emph{#1}}\hspace{\fill}$\square$\par}
\newcommand{\mcoro}[1]{%
\refstepcounter{snomer}
\par\textbf{Следствие \arabic{bnomer}.\arabic{snomer}. }{\emph{#1}}\par}
\newcommand{\corop}[2]{%
\refstepcounter{snomer}
\par\textbf{Следствие \arabic{bnomer}.\arabic{snomer}. }{\emph{#1}}
\par\textsc{Proof}. {#2}\hspace{\fill}$\square$\par}
\newcommand{\nota}[1]{%
\refstepcounter{snomer}
\par\medskip\textbf{Remark \arabic{bnomer}.\arabic{snomer}. }{#1}\par\medskip}
\newcommand{\propp}[2]{%
\refstepcounter{snomer}
\par\medskip\textbf{Proposition \arabic{bnomer}.\arabic{snomer}. }{\emph{#1}}
\par\textsc{Proof}. {#2}\hspace{\fill}$\square$\par\medskip}
\newcommand{\hypo}[1]{%
\refstepcounter{snomer}
\par\medskip\textbf{Conjecture \arabic{bnomer}.\arabic{snomer}. }{\emph{#1}}\par\medskip}

\newcommand{\Ind}[3]{%
\mathrm{Ind}_{#1}^{#2}{#3}}
\newcommand{\Res}[3]{%
\mathrm{Res}_{#1}^{#2}{#3}}
\newcommand{\epsi}{\epsilon}
\newcommand{\tri}{\triangleleft}
\newcommand{\Supp}[1]{%
\mathrm{Supp}(#1)}

\newcommand{\reg}{\mathrm{reg}}
\newcommand{\sreg}{\mathrm{sreg}}
\newcommand{\codim}{\mathrm{codim}\,}
\newcommand{\chara}{\mathrm{char}\,}
\newcommand{\rk}{\mathrm{rk}\,}
\newcommand{\chr}{\mathrm{ch}\,}
\newcommand{\id}{\mathrm{id}}
\newcommand{\Ad}{\mathrm{Ad}}
\newcommand{\col}{\mathrm{col}}
\newcommand{\row}{\mathrm{row}}
\newcommand{\low}{\mathrm{low}}
\newcommand{\pho}{\hphantom{\quad}\vphantom{\mid}}
\newcommand{\fho}[1]{\vphantom{\mid}\setbox0\hbox{00}\hbox to \wd0{\hss\ensuremath{#1}\hss}}
\newcommand{\wt}{\widetilde}
\newcommand{\wh}{\widehat}
\newcommand{\ad}[1]{\mathrm{ad}_{#1}}
\newcommand{\tr}{\mathrm{tr}\,}
\newcommand{\GL}{\mathrm{GL}}
\newcommand{\SL}{\mathrm{SL}}
\newcommand{\Sp}{\mathrm{Sp}}
\newcommand{\Mat}{\mathrm{Mat}}
\newcommand{\Stab}{\mathrm{Stab}}

\newcommand{\vfi}{\varphi}
\newcommand{\teta}{\vartheta}
\newcommand{\Fp}{\mathbb{F}}
\newcommand{\Rp}{\mathbb{R}}
\newcommand{\Zp}{\mathbb{Z}}
\newcommand{\Cp}{\mathbb{C}}
\newcommand{\ut}{\mathfrak{u}}
\newcommand{\at}{\mathfrak{a}}
\newcommand{\nt}{\mathfrak{n}}
\newcommand{\mt}{\mathfrak{m}}
\newcommand{\htt}{\mathfrak{h}}
\newcommand{\spt}{\mathfrak{sp}}
\newcommand{\rt}{\mathfrak{r}}
\newcommand{\rad}{\mathfrak{rad}}
\newcommand{\bt}{\mathfrak{b}}
\newcommand{\gt}{\mathfrak{g}}
\newcommand{\vt}{\mathfrak{v}}
\newcommand{\pt}{\mathfrak{p}}
\newcommand{\Xt}{\mathfrak{X}}
\newcommand{\Po}{\EuScript{P}}
\newcommand{\Uo}{\EuScript{U}}
\newcommand{\Fo}{\EuScript{F}}
\newcommand{\Do}{\EuScript{D}}
\newcommand{\Eo}{\EuScript{E}}
\newcommand{\Iu}{\mathcal{I}}
\newcommand{\Mo}{\mathcal{M}}
\newcommand{\Nu}{\mathcal{N}}
\newcommand{\Ro}{\mathcal{R}}
\newcommand{\Co}{\mathcal{C}}
\newcommand{\Lo}{\mathcal{L}}
\newcommand{\Ou}{\mathcal{O}}
\newcommand{\Au}{\mathcal{A}}
\newcommand{\Vu}{\mathcal{V}}
\newcommand{\Bu}{\mathcal{B}}
\newcommand{\Sy}{\mathcal{Z}}
\newcommand{\Sb}{\mathcal{F}}
\newcommand{\Gr}{\mathcal{G}}

\author{D.Y. Eliseev, M.V. Ignatyev\thanks{The second author was partially supported by RFBR grant no. 12--01--90805--mol\_rf\_nr.}
}
\date{\small Samara State University\\
Chair of algebra and geometry\\
\texttt{dmitriyelis@gmail.com}, \texttt{mihail.ignatev@gmail.com}}
\title{\Large{Kostant--Kumar polynomials and tangent cones\\ to Schubert varieties for involutions in $A_n$, $F_4$ and $G_2$}} \maketitle

\begin{flushright}
\textsl{To dear Professor Nikolai Vavilov\\with gratitude and admiration}
\end{flushright}

\sect{Introduction and the main result}

\sst Let $G$ be a complex reductive algebraic group, $T$ a maximal torus of $G$, $B$ a Borel subgroup of $G$ containing $T$, $\Phi$ the root system of $G$ w.r.t. $T$, $\Phi^+$ the set of positive roots w.r.t. $B$, $\Delta$ the set of fundamental roots, $W$ the Weyl group of $\Phi$ (see \cite{Bourbaki}, \cite{Humphreys} and \cite{Humpreys2} for basic facts about algebraic groups and root systems). Let $\Fo=G/B$ be the flag variety and $X_w\subseteq\Fo$ the Schubert subvariety corresponding to an element $w$ of the Weyl group $W$.

We denote by $\Ou=\Ou_{p,X_w}$ the local ring at the point $p=eB\in X_w$. Let $\mt$ be the maximal ideal of~$\Ou$. The sequence of ideals $$\Ou\supseteq\mt\supseteq\mt^2\supseteq\ldots$$ is a filtration. We define $\mathrm{gr}\,\Ou$ to be the graded algebra $$R=\mathrm{gr}\,\Ou=\bigoplus_{i\geq0}\mt^i/\mt^{i+1}.$$ By definition, the \emph{tangent cone} $C_w$ to the Schubert variety $X_w$ at the point $p$ is the spectrum of $R$: $C_w=\mathrm{Spec}\,R$. Clearly, $C_w$ is a subscheme of the tangent space $T_pX_w\subseteq T_p\Fo$. A natural problem in studying geometry of $X_w$ is to describe $C_w$ \cite[Chapter 7]{BilleyLakshmibai}.

Let $\gt$, $\bt$, $\htt$ be the Lie algebras of the groups $G$, $B$, $T$ respectively, $\htt^*$ the dual space of $\htt$. To each element $w\in W$ one can assign the polynomial $d_w\in S=\Cp[\htt]$ (see the next Subsection and \cite{KostantKumar1}, \cite{KostantKumar2}, \cite{Billey}, \cite[Section 7.1]{BilleyLakshmibai} for precise definitions). These polynomials are called Kostant--Kumar polynomials. They are the main tool in our study of tangent cones. In the paper \cite{Kumar}, S. Kumar showed that $d_w$ depends \emph{only} on $C_w$ (see the next Subsection for the details). In particular, to prove that $C_w$ and $C_{w'}$ do not coincide as subschemes of $T_p\Fo$, it is enough to check that $d_w\neq d_{w'}$.

In the paper \cite{EliseevPanov}, A.N. Panov and the first author computed tangent cones $C_w$ for all $w\in W$ in the case $G=\mathrm{SL}_n(\mathbb{C})$, $n\leq5$. Using this computations, Panov formulated several conjectures about the structure of tangent cones. In particular, he conjectured that if $C_w$ coincides with $C_{w'}$, then $w$ and $w'$ are conjugate in $W\cong S_n$. It turns out that this conjecture is false, $w=\begin{pmatrix}1&2&3&4&5&6\\3&5&4&6&2&1\end{pmatrix}$ and $w'=\begin{pmatrix}1&2&3&4&5&6\\3&6&5&2&1&4\end{pmatrix}$ give a counterexample (see also \cite{Bochkarev}). On the other hand, we have the following conjecture. \hypo{Suppose $w$\textup, $w'\in S_n$ are conjugate and $C_w=C_{w'}$. One can represent $w$ as a product of non-trivial disjoint cycles\textup: $w=c_1\ldots c_r$. Given $i$\textup, denote $F_i=\{j\mid c_i(j)\neq j\}\subseteq\{1,\ldots,n\}$. Let $H$ be the subgroup of $S_n$ consisting of all $g\in S_n$ satisfying $g(j)\in F_i$ for all $j\in F_i$\textup, $1\leq i\leq r$. Then there exists $g\in H$ such that $w'=gwg^{-1}$.}

This conjecture implies that if $w$ and $w'$ are involutions in $S_n$ (i.e., elements of order two) and $w\neq w'$, then $C_w\neq C_{w'}$. One can formulate the analogue of this conjecture for arbitrary $G$. (See papers \cite{Ignatyev1}, \cite{Ignatyev2} of the second author and Subsection~\ref{sst:coadjoint} for the interrelations between tangent cones to Schubert varieties $X_w$ associated with involutions and coadjoint orbits of the unipotent radical of the group~$B$.) To verify this conjecture, it is enough to check that the Kostant--Kumar polynomials of distinct involutions do not coincide. We prove this fact for $A_n$, $G_2$ and $F_4$. In the latter case, we use the computer algebra system \texttt{SAGE} \cite{SAGE}, see Subsection \ref{sst:F_4_G_2}. Precisely, we have the following result.
\mtheo{Assume that every irreducible component of the root system $\Phi$ is of type $A_n$\textup, $n\geq1$\textup, $F_4$ or $G_2$. Let $w$\textup, $w'$ be involutions in the Weyl group $W$ and $w\neq w'$. Then their Kostant--Kumar polynomials do not coincide\textup, i.e.\textup, $d_w\neq d_{w'}$. In particular\textup, the tangent cones $C_w$ and $C_{w'}$ do not coincide as subschemes of $T_p\Fo$.\label{mtheo}}

The structure of the paper is as follows. In Subsection~\ref{sst:definitions}, we give all necessary definitions and describe connections between tangent cones and Kostant--Kumar polynomials. In Subsections~\ref{sst:Bruhat}, \ref{sst:irred}, we recall some facts about the Bruhat order on the Weyl group and reduce the problem to the case of irreducible root systems (Proposition~\ref{post:irred}). Section~\ref{sect:proofs} contains the proof of the main theorem. In Subsections \ref{sst:A_n_1}--\ref{sst:A_n_2}, we prove it for $A_n$ by induction on $n$. In Subsection~\ref{sst:F_4_G_2}, we prove it for $F_4$ and~$G_2$. Section~\ref{sect:conclusion} contains some final remarks and conjectures.

A short announcement of our results was done in \cite{EliseevIgnatyev}.

\medskip\textsc{Acknowledgements}. We thank Professor Alexander N. Panov for useful discussions. A part of this work was carried out during the stay of the second author at Moscow State University. The second author thanks Professor Ernest B. Vinberg for his hospitality. Financial support from RFBR (grant no. 12--01--90805--mol\_rf\_nr) is gratefully acknowledged.

\sst\label{sst:definitions} Here we give precise definition of the Kostant--Kumar polynomial $d_w$, explain how to compute it in combinatorial terms and show that it depends only on the tangent cone $C_w$.

The torus $T$ acts on the Schubert variety $X_w$ by conjugation. Note that this action is the same as the action by left multiplication. The point $p$ is invariant under this action, so we have the structure of a $T$-module on the local ring $\Ou$. Clearly, the action of $T$ on $\Ou$ preserves the filtration by powers of the ideal $\mt$, hence we obtain the structure of a $T$-module on the algebra $R=\mathrm{gr}\,\Ou$. According to~\cite[Theorem~2.2]{Kumar}, $R$ can be decomposed into a direct sum of its finite-dimensional weight subspaces: $$R=\bigoplus_{\lambda\in\Xt(T)}R_{\lambda}.$$
Here $\Xt(T)\subseteq\htt^*$ is the character lattice of the torus $T$ and $R_{\lambda}=\{f\in R\mid t.f=\lambda(t)f\}$ is the weight subspace of weight $\lambda$. Let $\Lambda$ be the $\Zp$-module consisting of all (possibly infinite) $\Zp$-linear combinations of linearly independent elements $e^{\lambda}$, $\lambda\in\Xt(T)$. Then one can define the \emph{formal character} of $R$ to be an element of $\Lambda$ of the form $$\chr R=\sum_{\lambda\in\Xt(T)}m_{\lambda}e^{\lambda},$$ where $m_{\lambda}=\dim R_{\lambda}$.

Now, pick an element $a=\sum_{\lambda\in\Xt(T)}n_{\lambda}e^{\lambda}\in\Lambda$. Assume that there are finitely many $\lambda\in\Xt(T)$ such that $n_{\lambda}\neq0$. Given $k\geq0$, one can define the polynomial $$[a]_k=\sum_{\lambda\in\Xt(T)}n_{\lambda}\cdot\dfrac{\lambda^k}{k!}\in S=\Cp[\htt].$$ Denote $[a]=[a]_{k_0}$, where $k_0$ is minimal among all non-negative numbers $k$ such that $[a]_k\neq0$. For instance, if $a=1-e^{\lambda}$, then, $[a]_0=0$ and $[a]=[a]_1=-\lambda$ (here we denote $1=e^0$). Let $A$ be the submodule of $\Lambda$ consisting of all finite linear combinations. It is a commutative ring with respect to the multiplication $e^{\lambda}\cdot e^{\mu}=e^{\lambda+\mu}$. In fact, it is just the group ring of $\Xt(T)$. By $Q\subseteq\Lambda$, denote the field of fractions of the ring $A$. Note that to each element of $Q$ of the form $q=a/b$, $a$, $b\in A$, one can assign the element $$[q]=\dfrac{[a]}{[b]}\in\Cp(\htt)$$ of the field of rational functions on $\htt$.

There exists the involution $q\mapsto q^*$ on $Q$ defined by $$e^{\lambda}\mapsto (e^{\lambda})^*=e^{-\lambda}.$$ It turns out \cite[Theorem 2.2]{Kumar} that the character $\chr R$ belongs to $Q$, hence $(\chr R)^*\in Q$, too. Finally, we put $$c_w=[(\chr R)^*],\ d_w=(-1)^{l(w)}\cdot c_w\cdot\prod_{\alpha\in\Phi^+}\alpha.$$ Here $l(w)$ is the length of $w$ in the Weyl group $W$ with respect to the set of fundamental roots $\Delta$. Evidently, $c_w$ and $d_w$ belong to $\Cp(\htt)$; in fact, $d_w$ is a polynomial, i.e., belongs to $S=\Cp[\htt]$ (see \cite{KostantKumar2} and \cite[Theorem 7.2.6]{BilleyLakshmibai}). \defi{Let $w$ be an element of the Weyl group $W$. The polynomial $d_w\in S$ is called the \emph{Kostant--Kumar polynomial} associated with $w$.}

It follows from the definition that $c_w$ and $d_w$ depend only on the canonical structure of a $T$-module on the algebra $R$ of regular functions on the tangent cone $C_w$. Thus, to prove that the tangent cones corresponding to elements $w$, $w'$ of the Weyl group are distinct, it is enough to check that ${c_w\neq c_{w'}}$, or, equivalently, $d_w\neq d_{w'}$. On the other hand, there is a purely combinatorial description of Kostant--Kumar polynomials. To give this description, we need some more notation. Let $w$, $v$ be elements of~$W$. Fix a reduced expression of the element $w=s_{i_1}\ldots s_{i_l}$. (Here $\alpha_1,\ldots,\alpha_n\in\Delta$ are fundamental roots and $s_i$ is the simple reflection corresponding to $\alpha_i$.) Put
\begin{equation*}
c_{w,v}=(-1)^{l(w)}\cdot\sum\dfrac{1}{s_{i_1}^{\epsi_1}\alpha_{i_1}}\cdot\dfrac{1}{s_{i_1}^{\epsi_1}
s_{i_2}^{\epsi_2}\alpha_{i_2}}\cdot\ldots\cdot\dfrac{1}{s_{i_1}^{\epsi_1}\ldots s_{i_l}^{\epsi_l}\alpha_{i_l}},
\end{equation*}
where the sum is taken over all sequences $(\epsi_1,\ldots,\epsi_l)$ of zeroes and units such that $s_{i_1}^{\epsi_1}\ldots s_{i_l}^{\epsi_l}=v$. Actually, the element $c_{w,v}\in\Cp(\htt)$ depends only on $w$ and $v$, not on the choice of a reduced expression of $w$ \cite[Section 3]{Kumar}.

\exam{Let $\Phi=A_2$, so $W\cong S_3$. Put $w=s_1s_2s_1$. Let $\id$ be the identity element of the Weyl group. To compute $c_{w,\id}$, we should take the sum over two sequences, $(0,0,0)$ and $(1,0,1)$. Hence
\begin{equation*}
c_{w,\id}=(-1)^3\cdot\left(\dfrac{1}{\alpha_1\alpha_2\alpha_1}+\dfrac{1}{-\alpha_1(\alpha_1+\alpha_2)\alpha_1}\right)
=\dfrac{1}{\alpha_1\alpha_2(\alpha_1+\alpha_2)}.
\end{equation*}}

Now, put
\begin{equation}
d_{w,v}=\sum r(j_1)\ldots r(j_t)\in\Cp[\htt].\label{formula:d_w_v}
\end{equation}
Here $r(k)=s_{i_1}\ldots s_{i_{k-1}}\alpha_{i_k}$, $1\leq k\leq l$, and the sum is taken over all sequences $(j_1,\ldots,j_t)$, $t=l(v)$, such that $s_{i_{j_1}}\ldots s_{i_{j_t}}$ is a reduced expression of $v$ obtained from the expression $w=s_{i_1}\ldots s_{i_l}$ by deleting some simple reflections. Denote by $w_0$ the longest element of the Weyl group $W$. A remarkable fact is that \cite{KostantKumar2} \begin{equation}d_{vw_0,ww_0}=(-1)^{l(w)-l(v)}\cdot c_{w,v}\cdot\prod_{\alpha\in\Phi^+}\alpha.\label{formula:d_and_c}\end{equation} In particular, $d_{w,v}$ does not depend on the choice of a reduced expression of $w$. Further, $c_w=c_{w,\id}$ (and so $d_w=d_{w_0,ww_0}$), hence to prove that tangent cones do not coincide, we need only combinatorics of the Weyl group.

\medskip At the rest of the Subsection, we present an original definition of elements $c_{w,v}$ using so-called nil-Hecke rings (see \cite{Kumar} and~\cite[Section 7.1]{BilleyLakshmibai}). (This definition is needed for us in the case $A_n$.) Denote by $Q_W$ the vector space over $\Cp(\htt)$ with basis $\{\delta_w,\ w\in W\}$. It is a ring with respect to the multiplication $$f\delta_v\cdot g\delta_w=fv(g)\delta_{vw}.$$
This ring is called the \emph{nil-Hecke ring}. To each $i$ from 1 to $n$ put $$x_i=\alpha_i^{-1}(\delta_{s_i}-\delta_{\id}).$$
Let $w\in W$ and $w=s_{i_1}\ldots s_{i_l}$ be a reduced expression of $W$. Then the element $$x_w=x_{i_1}\ldots x_{i_l}$$ does not depend on the choice of a reduced expression \cite[Proposition 2.1]{KostantKumar1}.

Furthermore, it turns out that $\{x_w,\ w\in W\}$ is a $\Cp(\htt)$-basis of~$Q_W$ \cite[Proposition 2.2]{KostantKumar1}, and
\begin{equation*}
\begin{split}
&x_w=\sum\nolimits_{v\in W}c_{w,v}\delta_v,\\
&\delta_w=\sum\nolimits_{v\in W}d_{w,v}x_v.
\end{split}
\end{equation*}
Actually, if $w,v\in W$, then
\begin{equation}
\begin{split}
&\text{a) }x_v\cdot x_w=\begin{cases}x_{vw},&\text{ if }l(vw)=l(v)+l(w),\\
0,&\text{ otherwise},
\end{cases}\\
&\text{b) }c_{w,v}=-v(\alpha_i)^{-1}(c_{ws_i,v}+c_{ws_i,vs_i}),\text{ if }l(ws_i)=l(w)-1.\\
\end{split}\label{formula_x_v_x_w}
\end{equation}
(The group $W$ naturally acts on $\Cp(\htt)$ by automorphisms.) The first property is proved in \cite[Proposition 2.2]{KostantKumar1} and the second property follows immediately from the first one and the definitions (see also the proof of \cite[Corollary 3.2]{Kumar}).

\sst\label{sst:Bruhat} In this Subsection, we briefly recall some facts about the Bruhat order on the Weyl group needed for the sequel. We say that $v$ is less or equal to $w$ with respect to the Bruhat order, written $v\leq w$, if some reduced expression for $v$ is a subword of some reduced expression for $w$. It is well-known that this order plays the crucial role in many geometric aspects of theory of algebraic groups. For instance, the Bruhat order encodes the incidences among Schubert varieties, i.e., $X_v$ is contained in the closure of $X_w$ if and only if $v\leq w$.

It turns out that $c_{w,v}$  is non-zero if and only if $v\leq w$ \cite[Corollary 3.2]{Kumar}. For example, $c_w=c_{w,\id}$ iz non-zero for \emph{all} $w$, because $\id$ is the smallest element of $W$ with respect to the Bruhat order. Note that (see \cite{Dyer} and \cite[Theorem 7.1.11]{BilleyLakshmibai}) given $v,w\in W$, there exists $g_{w,v}\in S=\Cp[\htt]$ such that
\begin{equation}
c_{w,v}=g_{w,v}\cdot\prod_{\alpha>0,\ s_{\alpha}v\leq w}\alpha.\label{formula:dyer}
\end{equation}

In Subsections \ref{sst:A_n_1}, \ref{sst:A_n_2}, we will use the Bruhat order on the symmetric group. In this case, there exists a nice description of the Bruhat order. Namely, given $w\in S_n$, denote by $\dot w$ the $n\times n$ matrix of the form
\begin{equation*}
(\dot w)_{i,j}=\begin{cases}1,&\text{ if }w(j)=i,\\
0&\text{ otherwise}.
\end{cases}
\end{equation*}
It it called $0$--1 matrix, permutation matrix or rook placement for $w$. Define the matrix $R_w$ by putting its $(i,j)$th element to be equal to the rank of the lower left $(n-i+1)\times j$ submatrix of $\dot w$. In other words, $(R_w)_{i,j}$ is just the number or rooks located non-strictly to the South-West from $(i,j)$.

\exam{Let $n=7$, $w=\begin{pmatrix}1&2&3&4&5&6&7\\1&6&2&3&4&5&7\end{pmatrix}$. Here we draw the matrices $\dot w$ and $R_w$ (rooks are marked by $\otimes$):
\begin{equation*}\dot w=
\mymatrix{
\otimes& \pho& \pho& \pho& \pho& \pho& \pho\\
\pho& \pho& \otimes& \pho& \pho& \pho& \pho\\
\pho& \pho& \pho& \otimes& \pho& \pho& \pho\\
\pho& \pho& \pho& \pho& \otimes& \pho& \pho\\
\pho& \pho& \pho& \pho& \pho& \otimes& \pho\\
\pho& \otimes& \pho& \pho& \pho& \pho& \pho\\
\pho& \pho& \pho& \pho& \pho& \pho& \otimes\\
}\ ,\ R_w=
\mymatrix{
\fho1& \fho2& \fho3& \fho4& \fho5& \fho6& \fho7\\
\fho0& \fho1& \fho2& \fho3& \fho4& \fho5& \fho6\\
\fho0& \fho1& \fho1& \fho2& \fho3& \fho4& \fho5\\
\fho0& \fho1& \fho1& \fho1& \fho2& \fho3& \fho4\\
\fho0& \fho1& \fho1& \fho1& \fho1& \fho2& \fho3\\
\fho0& \fho1& \fho1& \fho1& \fho1& \fho1& \fho2\\
\fho0& \fho0& \fho0& \fho0& \fho0& \fho0& \fho1\\
}\ .
\end{equation*}\label{exam:Bruhat}}

Let $X$ and $Y$ be matrices with integer entries. We say that $X\leq Y$ if $X_{i,j}\leq Y_{i,j}$ for all $i,j$. It turns out that if $v$, $w\in S_n$, then
\begin{equation}
v\leq w \text{ if and only if }R_v\leq R_w\label{formula:Bruhat}
\end{equation}(see, e.g., \cite[Theorem 1.6.4]{Incitti}).

\sst\label{sst:irred} Here we explain why it is enough to prove Theorem~\ref{mtheo} for irreducible root systems. It follows immediately from the next Proposition. Suppose~$\Phi$ is a union of its subsystems $\Phi_1$ and $\Phi_2$ contained in mutually orthogonal subspaces. Let $W_1$, $W_2$ be the Weyl groups of $\Phi_1$, $\Phi_2$ respectively, so $W=W_1\times W_2$. Denote $\Delta_1=\Delta\cap\Phi_1=\{\alpha_1,\ldots,\alpha_r\}$ and $\Delta_2=\Delta\cap\Phi_2=\{\beta_1,\ldots,\beta_s\}$, then $$S=\Cp[\htt]\cong\Cp[\alpha_1,\ldots,\alpha_r,\beta_1,\ldots,\beta_s].$$

Given $v\in W_1$, denote by $d_v^1$ its Kostant--Kumar polynomial. We can consider $d_v^1$ as an element of~$S$ depending only on $\alpha_1,\ldots,\alpha_r$. We define $c_v^1\in\Cp(\htt)$ by the similar way. Given $v\in W_2$, we define $d_v^2\in\Cp[\htt]$ and $c_v^2\in\Cp(\htt)$; they depend only on $\beta_1,\ldots,\beta_s$.

\postp{Let $w\in W$, $w_1\in W_1$, $w_2\in W_2$ and $w=w_1w_2$. \label{post:irred}Then $$d_w=d_{w_1}^1d_{w_2}^2,\ c_w=c_{w_1}^1c_{w_2}^2.$$}{By $s_i$ (resp. $r_j$), we denote the simple reflection corresponding to a simple root $\alpha_i$ (resp. $\beta_j$). Let $l_i$ be the length function on $W_i$ with respect to $\Delta_i$, $i=1,2$. It is well-known that $l_i(v)=l(v)$ for all $v\in W_i$. Hence if
\begin{equation*}
w_1=s_{i_1}\ldots s_{i_p},\ w_2=r_{j_1}\ldots r_{j_q}
\end{equation*}
are reduced expressions for $w_i$ in $W_i$, then they are reduced expressions for $w_i$ in $W$. Moreover, $$l(w)=l(w_1)+l(w_2)=l_1(w_1)+l_2(w_2).$$ This means that
$$w=s_{i_1}\ldots s_{i_p}r_{j_1}\ldots r_{j_q}$$
is a reduced expression for $w$ in $W$.\newpage

It follows from $W=W_1\times W_2$ that $$s_{i_1}^{\epsi_1}\ldots s_{i_p}^{\epsi_p}r_{j_1}^{\delta_1}\ldots r_{j_q}^{\delta_q}=\id,$$ $\epsi_i,\delta_j\in\{0,1\}$, is equivalent to $$s_{i_1}^{\epsi_1}\ldots s_{i_p}^{\epsi_p}=r_{j_1}^{\delta_1}\ldots r_{j_q}^{\delta_q}=\id.$$ Since all $s_i$'s (resp. $r_j$'s) act identically on $\Phi_2$ (resp. on $\Phi_1$), we obtain
\begin{equation*}
\begin{split}
c_w=(-1)^{l_1(w_1)+l_2(w_2)}&\cdot\sum\left(\dfrac{1}{s_{i_1}^{\epsi_1}\alpha_{i_1}}\cdot\dfrac{1}{s_{i_1}^{\epsi_1}
s_{i_2}^{\epsi_2}\alpha_{i_2}}\cdot\ldots\cdot\dfrac{1}{s_{i_1}^{\epsi_1}\ldots s_{i_p}^{\epsi_p}\alpha_{i_p}}\vphantom{\dfrac{1}{r_{j_1}^{\delta_1}\ldots r_{j_q}^{\delta_q}\beta_{j_q}}}\right.\\
&\times\left.\dfrac{1}{r_{j_1}^{\delta_1}\beta_{j_1}}\cdot\dfrac{1}{r_{j_1}^{\delta_1}
r_{j_2}^{\delta_2}\beta_{j_2}}\cdot\ldots\cdot\dfrac{1}{r_{j_1}^{\delta_1}\ldots r_{j_q}^{\delta_q}\beta_{j_q}}\right)=c_{w_1}^1c_{w_2}^2.
\end{split}
\end{equation*}
The second equality is proved. The first equality follows immediately from the second one and the obvious fact that $\Phi^+=\Phi_1^+\cup\Phi_2^+$.}

Now, to prove the main Theorem, it suffice to check it for irreducible root systems of types $A_n$, $F_4$ and $G_2$, because $\Cp[\htt]$ is a unique factorization domain.

\sect{Proofs}\label{sect:proofs}

\sst\label{sst:A_n_1} In this Subsection and in the next Subsection, we will prove the main result for the case $\Phi=A_{n-1}$, $n\geq2$. As usual, we identify $\Phi^+$ with the subset of the Euclidean space $\Rp^n$ of the form $$\{\epsi_j-\epsi_i,\ 1\leq j<i\leq n\}$$ ($\{\epsi_i\}_{i=1}^n$ is the standard basis). In this case, $W$ is isomorphic to $S_n$, the symmetric group on $n$ letters, and a transposition $(i,j)$ is just a reflection $s_{\epsi_j-\epsi_i}$. Here $\alpha_1=\epsi_1-\epsi_2$, $\ldots$, $\alpha_{n-1}=\epsi_{n-1}-\epsi_n$ are fundamental roots.

We will consider not all elements of $W$, but only \emph{involutions}, i.e., elements of order two. We put $$I_n=I(W)=\{\sigma\in W\mid\sigma^2=\id\}.$$ Each involution $\sigma$ can be uniquely presented as a product of disjoint 2-cycles $\sigma=(i_1,j_1)\ldots(i_l,j_l)$, $i_k>j_k$, $j_1<\ldots<j_l$. \defi{The \emph{support} of an involution $\sigma\in I_n$ is the subset of $\Phi^+$ of the form $$\Supp{\sigma}=\{\epsi_{j_1}-\epsi_{i_1},\ldots,\epsi_{j_l}-\epsi_{i_l}\}.$$ Note that it consists of pairwise orthogonal roots. In other words, the support of $\sigma$ is the unique orthogonal subset of $\Phi^+$ such that $$\sigma=\prod_{\alpha\in\Supp{\sigma}}s_{\alpha}.$$ (Here reflections are taken in any fixed order: since the support is an orthogonal subset, they commute.)}\exam{If $n=7$ and $\sigma=\begin{pmatrix}1&2&3&4&5&6&7\\5&7&6&4&1&3&2\end{pmatrix}=(5,1)(7,2)(6,3)$, then $$\Supp{\sigma}=\{\epsi_1-\epsi_5,\ \epsi_2-\epsi_7,\ \epsi_3-\epsi_6\}.$$}

Note also that there is a quite simple description of the Bruhat order on involutions in $S_n$. Namely, let $w\in I_n$. Let $R_w$ be the matrix defined in Subsection~\ref{sst:Bruhat}, and $R_w^*$ its strictly lower-triangular part, i.e.,
\begin{equation*}
(R_w^*)_{i,j}=\begin{cases}(R_w)_{i,j},&\text{ if }i>j,\\
0&\text{ otherwise.}\\
\end{cases}
\end{equation*}
By \cite[Theorem 1.10]{Ignatyev2}, if $v,w\in I_n$, then
\begin{equation}
v\leq w \text{ if and only if }R_v^*\leq R_w^*.\label{formula:Bruhat_inv}
\end{equation}

We will prove Theorem~\ref{mtheo} by induction on $n$ (for $n=2$, there is nothing to prove). Denote by $\wt W=\wt S_{n-1}$ the subgroup of $W$ consisting of all permutations $w$ such that $w(1)=1$; clearly, $\wt W\cong S_{n-1}$. Let $\wt I_{n-1}=I(\wt W)$ be the set of involutions in $\wt W$. Given $w\in\wt W$, we denote by $\wt d_{w}$ its Kostant--Kumar polynomial. One can identify $\Cp[\htt]$ with $\Cp[\alpha_1,\ldots,\alpha_{n-1}]$, then $\wt d_w$ belongs to $\Cp[\htt]$ and does not depend on $\alpha_1$. Similarly, we define $\wt c_w\in\Cp(\htt)$ and $\wt d_{w,v}$, $\wt c_{w,v}$ for all $w,v\in W$. By the inductive assumption, $\wt d_w\neq\wt d_v$ and $\wt c_w\neq\wt c_w$ for all distinct involutions $w,v\in\wt I_{n-1}$.

We need some more notation. For any $\alpha=\epsi_j-\epsi_i\in\Phi^+$, define $\row(\alpha)=i$, $\col(\alpha)=j$. For any $k$ from 1 to $n$, put
\begin{equation*}
\begin{split}
\Ro_k&=\{\alpha\in\Phi^+\mid\row(\alpha)=k\},\\
\Co_k&=\{\alpha\in\Phi^+\mid\col(\alpha)=k\}.\\
\end{split}
\end{equation*}
The set $\Ro_k$ (resp. $\Co_k$) is called the $k$th \emph{row} (resp. the $k$th \emph{column}). We have $$\wt I_{n-1}=\{\sigma\in I_n\mid\Supp{\sigma}\cap\Co_1=\varnothing\}.$$ Furthermore, for any $k$ and any involution $\sigma\in I_n$,
\begin{equation*}
|\Supp{\sigma}\cap(\Ro_k\cup\Co_k)|\leq1.
\end{equation*}

\nota{There is a natural order on the root system $\Phi$. By definition, $\alpha\leq\beta$ means that $\beta-\alpha$ is a sum of positive roots. In other words, $\alpha=\epsi_j-\epsi_i\leq\beta=\epsi_s-\epsi_r$ if and only if $s\leq j$ and~$i\leq r$. Using~(\ref{formula:Bruhat_inv}), one can easily check that if $w$ is an involution and $\alpha=\epsi_j-\epsi_i$ is a positive root, then $s_{\alpha}\leq w$ if and only if $\alpha\leq\beta$ for some positive root $\beta=\epsi_s-\epsi_r\in\Supp{\sigma}$. Indeed, suppose the latter condition holds. Then
\begin{equation*}
(R_{s_{\alpha}}^*)_{k,l}=\begin{cases}1,&\text{if }j\leq k<l\leq i,\\
0&\text{otherwise,}
\end{cases}
\end{equation*}
and $(R_w^*)_{k,l}\geq1$ for all $s\leq k<l\leq r$, so $s_{\alpha}\leq w$. At the contrary, if this condition does not hold, then $$(R_{s_{\alpha}}^*)_{i,j}=1>0=(R_w^*)_{i,j},$$
so $s_{\alpha}\not\leq w$. In particular, if $\Supp{\sigma}\cap\Co_1=\{\beta\}$, where $\beta=\epsi_1-\epsi_i$, and $\alpha=\epsi_1-\epsi_k\in\Co_1$, then $s_{\alpha}\leq\sigma$ if and only if $\alpha\leq\beta$, i.e., $k\leq i$.\label{nota:order_roots}}

Now we will prove two important Lemmas. \lemmp{Let $w\in\wt I_{n-1}$. Then $d_w=\wt d_w\cdot\prod_{\alpha\in\Co_1}\alpha$.\label{lemm:A_n_C_1_ne_C_1}}{Since $\wt W$ is a parabolic subgroup of $W$, any reduced expression for $w$ in $\wt W$ is a reduced expression for $w$ in $W$. This implies $\wt c_w=c_w$. The result follows.}\newpage

\lemmp{Let $w\in I_n$. Suppose $\Supp{w}\cap\Co_1=\{\beta\}$ and $$c_w=A/B,\ A,B\in\Cp[\htt],\ (A,B)=1,$$ i.e.\textup{,} $A$ and $B$ are coprime. Then $\beta$ divides $B$ in the polynomial ring $\Cp[\htt]$.\label{lemm:A_n_C_1}}{Suppose $\beta=\epsi_1-\epsi_j$. Put
$$u=s_{j-1}\ldots s_1=(j,j-1)\ldots(2,1)=\begin{pmatrix}1&2&3&\ldots&j-1&j&j+1&\ldots&n\\j&1&2&\ldots&j-2&j-1&j+1&\ldots&n\end{pmatrix}.$$
Denote $v=u^{-1}w$, so $w=uv$. Clearly, $v(1)=u^{-1}(w(1))=u^{-1}(j)=1$, so $v\in\wt W$. Further,
\begin{equation*}
\begin{split}
&u(\alpha_i)=u(\epsi_i-\epsi_{i+1})=\epsi_{i-1}-\epsi_i>0\text{ for all }i\text{ from 2 to }j-1,\\
&u(\alpha_j)=u(\epsi_j-\epsi_{j+1})=\epsi_{j-1}-\epsi_{j+1}>0,\\
&u(\alpha_i)=u(\epsi_i-\epsi_{i+1})=\epsi_i-\epsi_{i+1}=\alpha_i>0\text{ for all }i\text{ from $j+1$ to }n-1.
\end{split}
\end{equation*}
By the way, $u(\alpha_i)>0$ if $i\geq2$. This is equivalent to $l(us_i)=l(u)+1$. According to \cite[Pro\-po\-si\-tion~1.10]{Humpreys2}, $l(w)=l(u)+l(v)$.

Using~(\ref{formula_x_v_x_w}a), we obtain
\begin{equation*}
\begin{split}
x_w&=\sum_{s\in W}c_{w,s}\delta_s=x_ux_v=\sum_{g,h\in W}c_{u,g}\delta_g\cdot c_{v,h}\delta_h\\
&=\sum_{g,h\in W}c_{u,g}g(c_{v,h})\delta_{gh}=\sum_{s\in W}\left(\sum_{g\in W}c_{u,g}g(c_{v,g^{-1}s})\right)\delta_s.
\end{split}
\end{equation*}
Thus, for any $s\in W$, the coefficient of $\delta_s$ is equal to
$$c_{w,s}=\sum_{g\in W}c_{u,g}g(c_{v,g^{-1}s}),$$
in particular,
$$c_w=c_{w,\id}=\sum_{g\in W}c_{u,g}g(c_{v,g^{-1}}).$$
Moreover, since $c_{p,q}\neq0$ if and only if $p\geq q$, the sum in the right hand side is taken over permutations $g$ such that $u\geq g$ and $v\geq g^{-1}$. Denote the set of such permutations by $U$. Note that $g\in U$ implies that $g$ is obtained from $u=s_{j-1}\ldots s_1$ by deleting $s_1$ and, possibly, some other simple reflections. (If $s_1$ is not deleted, then the condition $v\geq g^{-1}$ does not hold.) Hence
$$c_w=c_{w,\id}=\sum_{g\in U}c_{u,g}g(c_{v,g^{-1}}).$$

Using (\ref{formula_x_v_x_w}b) and the fact that $l(us_1)=l(u)-1$, we obtain
$$c_{u,g}=-g(\alpha_1)^{-1}(c_{us_1},g+c_{us_1,gs_1})=-g(\alpha_1)^{-1}c_{us_1,g},$$
because $us_1\not\geq gs_1$ and so $c_{us_1,gs_1}=0$. Thus,
$$c_w=-\sum_{g\in U}\dfrac{c_{us_1,g}g(c_{v,g^{-1}})}{g\alpha_1}.$$
It is easy to check that there is most one $g$ such that $g\alpha_1=\beta$ and $g\in U$, namely, $g_0=us_1=s_{j-1}\ldots s_2$. Clearly, $g_0\alpha_1=\beta$. Assume for a moment that $g_0$ belongs to $U$, i.e., $v\geq g_0^{-1}$.

Then
\begin{equation}
c_w=-\dfrac{c_{us_1,g_0}g_0(c_{v,g_0^{-1}})}{\beta}-\sum_{g\in U,\ g\neq g_0}\dfrac{c_{us_1,g}g(c_{v,g^{-1}})}{g\alpha_1}.\label{formula:c_w_c_v}
\end{equation} By $S'$ (resp. $Q'$) denote the subalgebra of $S=\Cp[\htt]$ (resp. the subfield of $\Cp(\htt)$) generated by\linebreak $\alpha_2,\ldots,\alpha_{n-1}$, then $c_{v,g_0^{-1}}\in Q'$. Since $g(1)=1$, $g(c_{v,g_0^{-1}})\in Q'$, too. In particular, if $g(c_{v,g_0^{-1}})=G_1/G_2$ and $G_1$, $G_2\in S'$ are coprime, then $\beta$ does not divide $G_1$. On the other hand,
$$c_{us_1,g_0}=c_{us_1,us_1}=\pm\dfrac{1}{s_{j-1}\alpha_{j-1}}\cdot\dfrac{1}{s_{j-1}s_{j-2}\alpha_{j-2}}\cdot\ldots\cdot\dfrac{1}{s_{j-1}\ldots s_2\alpha_2},$$
because $us_1=s_{j-1}\ldots s_2$. We conclude that the first summand in the sum above has the form $P/\beta Q$ for some coprime $P,Q\in\Cp[\htt]$ such that $P$ is \emph{non-zero}.

Similarly, if $g\in U$ and $g\neq g_0$, then $g(c_{v,g^{-1}})\in Q'$. At the same time, $$c_{us_1,g}=\pm\dfrac{1}{s_{l_1}\alpha_{l_1}}\cdot\dfrac{1}{s_{l_1}s_{l_2}\alpha_{l_2}}\cdot\ldots\cdot\dfrac{1}{s_{l_1}\ldots s_{l_k}\alpha_{l_k}},$$
where $g=s_{l_1}\ldots s_{l_k}$ for certain $j-1\geq l_1>l_2>\ldots>l_k\geq2$. We see that if the latter sum in~\ref{formula:c_w_c_v} is equal to $C/D$, where $C,D\in\Cp[\htt]$ are coprime, then $\beta$ does not divide $D$. Thus,
$$c_w=\dfrac{C}{D}+\dfrac{P}{\beta Q}=\dfrac{\beta CQ+PD}{\beta DQ}.$$ Here $\beta$ divides neither $P$ nor $D$, hence $\beta$ does not divide the numerator. Thus, $\beta$ divides the denominator of $c_w$, as required.

Thus, to conclude the proof, we must show that $g_0\in U$, i.e., $v\geq g_0^{-1}$, or, equivalently, $v^{-1}\geq g_0$. To do this, note that
\begin{equation*}
(R_{g_0})_{p,q}=\begin{cases}p-q+1,&\text{ if }p\leq q,\\
1,&\text{ if }2\leq q<p\leq j,\\
0,&\text{ otherwise}.
\end{cases}
\end{equation*}
(In fact, Example~\ref{exam:Bruhat} deals with $g_0$ for $n=7$, $j=6$.) At the same time,
\begin{equation*}
v^{-1}=w^{-1}u=wu=\begin{pmatrix}1&\ldots&j&\ldots\\j&\ldots&1&\ldots\end{pmatrix}
\cdot\begin{pmatrix}1&2&\ldots\\j&1&\ldots\end{pmatrix}=\begin{pmatrix}1&2&\ldots\\1&j&\ldots\end{pmatrix},\\
\end{equation*}
so if $2\leq q<p\leq j$, then $(R_{v^{-1}})_{p,q}\geq1$. But if $1\leq p\leq q\leq n$, then $(R_{g_0})_{p,q}=(R_{\id})_{p,q}$. Since $\id$ is the smallest element of $W$ with respect to the Bruhat order, (\ref{formula:Bruhat}) shows that $v^{-1}\geq g_0$. The proof is complete.}

\sst\label{sst:A_n_2} Now, we can prove the main theorem for $A_n$. The prove follows immediately from Pro\-po\-si\-tions~\ref{prop:A_n_C_1_ne_C_1}, \ref{prop:A_n_C_1_C_1_neq} and~\ref{prop:A_n_C_1_C_1_eq}. Recall the notation from the previous Subsection. \propp{Let $\sigma,\tau\in I_n$ be involutions. Suppose $\Supp{\sigma}\cap\Co_1\neq\varnothing$ and $\Supp{\tau}\cap\Co_1=\varnothing$, then $d_{\sigma}\neq d_{\tau}$.\label{prop:A_n_C_1_ne_C_1}}{Suppose $\Supp{\sigma}\cap\Co_1=\{\beta\}$. Lemma~\ref{lemm:A_n_C_1_ne_C_1} shows that $\beta$ divides~$d_{\tau}$ in the polynomial ring~$\Cp[\htt]$. On the other hand, Lemma~\ref{lemm:A_n_C_1} claims that there exist coprime $A,B\in\Cp[\htt]$ such that $c_{\sigma}=A/B$, $\beta$ divides $B$ and does not divide $A$. Hence
$$d_{\sigma}=\pm c_{\sigma}\cdot\prod_{\alpha>0}\alpha=\pm\prod_{\alpha>0}\alpha\cdot A/B,$$ so $\beta$ does not divide $d_{\sigma}$. We conclude that $d_{\tau}\neq d_{\sigma}$, as required. Note that we did not use induction in this proof.}

\propp{Let $\sigma,\tau\in I_n$ be involutions. Suppose $\Supp{\sigma}\cap\Co_1=\{\beta\}$, $\Supp{\tau}\cap\Co_1=\{\gamma\}$ and $\beta\neq\gamma$, then $d_{\sigma}\neq d_{\tau}$.\label{prop:A_n_C_1_C_1_neq}}{Assume without loss of generality that $\beta>\gamma$, i.e., if $\beta=\epsi_1-\epsi_i$, $\gamma=\epsi_1-\epsi_s$, then $i>s$ (see Remark~\ref{nota:order_roots}). This Remark also shows that $s_{\beta}\not\leq\tau$. By formula~(\ref{formula:dyer}), there exists $g=g_{\tau,\id}\in\Cp[\htt]$ such that
\begin{equation*}
d_\tau=\pm c_{\tau}\cdot\prod_{\alpha>0}\alpha=\pm g\cdot\prod_{\alpha>0,\ s_{\alpha}\not\leq\tau}\alpha,
\end{equation*}
so $\beta$ divides $d_{\tau}$, because $\beta$ is involved in the latter product. As in the previous Proposition, using Lemma~\ref{lemm:A_n_C_1}, we obtain that $\beta$ does not divide $d_{\sigma}$. Thus, $d_{\tau}\neq d_{\sigma}$. Note that we did not use induction in this proof.}

\propp{Let $\sigma,\tau\in I_n$ be distinct involutions. Suppose $\Supp{\sigma}\cap\Co_1=\{\beta\}=\Supp{\tau}\cap\Co_1$, then $d_{\sigma}\neq d_{\tau}$.\label{prop:A_n_C_1_C_1_eq}}{Suppose $\beta=\epsi_1-\epsi_j$. Consider an involution $w\in I_n$ such that $\Supp{w}\cap\Co_1=\{\beta\}$. As in the proof of Lemma~\ref{lemm:A_n_C_1}, put $w=uv$, where $u=s_{j-1}\ldots s_1$ and $v=u^{-1}w\in\wt W$. Recall from (\ref{formula:c_w_c_v}) that
\begin{equation*}
c_w=-\dfrac{c_{us_1,g_0}g_0(c_{v,g_0^{-1}})}{\beta}-\sum_{g\in U,\ g\neq g_0}\dfrac{c_{us_1,g}g(c_{v,g^{-1}})}{g\alpha_1},
\end{equation*}
where $U=\{g\in W\mid g\leq u,\ g^{-1}\leq v\}$ and $g_0=us_1\in U$.

Now, denote $w'=s_{j-1}ws_{j-1}\in I_n$. Assume that $j>2$, then $\Supp{w'}\cap\Co_1=\beta'=\epsi_1-\epsi_{j-1}$. As above, put $w'=u'v'$, where $u'=s_{j-2}\ldots s_1$ and $v'\in\wt W$, then
\begin{equation*}
c_{w'}=-\dfrac{c_{u's_1,h_0}h_0(c_{v',h_0^{-1}})}{\beta}-\sum_{h\in U',\ h\neq h_0}\dfrac{c_{u's_1,h}h(c_{v,h^{-1}})}{h\alpha_1},
\end{equation*}
where $U'=\{h\in W\mid h\leq u',\ h^{-1}\leq v'\}$ and $h_0=u's_1\in U$.

Our goal is to compare $c_{v,g_0^{-1}}$ with $c_{v',h_0^{-1}}$. Note that $u'=s_{j-1}u$, $v'=vs_{j-1}$ and $h_0=s_{j-1}g_0$. Recall that $$u=s_{j-1}\ldots s_1=(j,j-1)\ldots(2,1)=\begin{pmatrix}1&2&3&\ldots&j-1&j&j+1&\ldots&n\\j&1&2&\ldots&j-2&j-1&j+1&\ldots&n\end{pmatrix},$$
hence
$$v\alpha_{j-1}=u^{-1}w(\epsi_{j-1}-\epsi_j)=u^{-1}(\epsi_x-\epsi_1)=\epsi_y-\epsi_2,$$
where $x=w(j-1)$, $y=u^{-1}(x)=v(j-1)$. If $y=1$, then $u^{-1}(x)=1$, so $x=j$, but $w(j-1)\neq j$, a~contradiction, Hence $y>2$, so $v\alpha_{j-1}<0$. This means that $l(vs_{j-1})=l(v)-1$. Formula (\ref{formula_x_v_x_w}b) implies that
$$c_{v,g_0^{-1}}=\dfrac{c_{vs_{j-1},g_0^{-1}}+c_{vs_{j-1},g_0^{-1}s_{j-1}}}{-g_0^{-1}\alpha_{j-1}}.$$

We see that $$g_0=s_{j-1}\ldots s_2=\begin{pmatrix}1&2&3&\ldots&j-1&j&j+1&\ldots&n\\1&j&2&\ldots&j-2&j-1&j+1&\ldots&n\end{pmatrix},$$
hence $(R_{g_0})_{j,2}=1$. On the other hand,
$$(vs_{j-1})^{-1}=s_{j-1}wu=\begin{pmatrix}1&2&\ldots\\1&j-1&\ldots\end{pmatrix},$$
hence $(R_{(vs_{j-1})^{-1}})_{j,2}=0$. Formula (\ref{formula:Bruhat}) claims that $(vs_{j-1})^{-1}\not\geq g_0$, or, equivalently, $vs_{j-1}\not\geq g_0^{-1}$. We obtain $c_{vs_{j-1},g_0^{-1}}=0$, so
$$c_{v,g_0^{-1}}=\dfrac{c_{vs_{j-1},g_0^{-1}s_{j-1}}}{-g_0^{-1}\alpha_{j-1}}=\dfrac{c_{v',h_{0}^{-1}}}{\epsi_2-\epsi_j}.$$

If $j-1>2$, then we repeat this procedure with $w'$ in place of $w$, etc. In a finite number of steps we will obtain $w=aw_1a^{-1}$, where $a=s_2s_3\ldots s_{j-1}$. Here $w_1\in I_n$ and $\Supp{w_1}\cap\Co_1=\{\alpha_1\}=\{\epsi_1-\epsi_2\}$. We denote $w_1=u_1v_1$, where $u_1=s_1$ and $v_1\in\wt W$ is an \emph{involution}, i.e., $v_1\in\wt I_{n-1}$. It follows from the above that $c_{v,g_0^{-1}}=fc_{v_1,\id}$, where $$f=\dfrac{1}{(\epsi_2-\epsi_j)\cdot(\epsi_2-\epsi_{j-1})\cdot\ldots\cdot(\epsi_2-\epsi_3)}$$ \emph{does not depend} on $w$.

Now, consider the involutions $\sigma$ and $\tau$. Put $\sigma=uv_{\sigma}$ and $\tau=uv_{\tau}$, as above. Since $\sigma\neq\tau$, $\sigma_1\neq\tau_1$, too, where $\sigma_1=a\sigma a^{-1}$, $\tau_1=a\tau a^{-1}$. Hence $v_{\sigma}^1=s_1\sigma_1\neq v_{\tau}^1=s_1\tau_1$. By the inductive assumption applied to $v_{\sigma}^1$, $v_{\tau}^1\in\wt I_{n-1}$, one has $\wt c_{v_{\sigma}^1,\id}\neq\wt c_{v_{\tau}^1,\id}$. Lemma~\ref{lemm:A_n_C_1_ne_C_1} says that $c_{v_{\sigma}^1,\id}=\wt c_{v_{\sigma}^1,\id}\neq\wt c_{v_{\tau}^1,\id}=c_{v_{\tau}^1,\id}$, and, consequently, $$c_{v_{\sigma},g_0^{-1}}=fc_{v_{\sigma}^1,\id}\neq fc_{v_{\tau}^1,\id}=c_{v_{\tau},g_0^{-1}},$$ hence $g_0(c_{v_{\sigma},g_0^{-1}})\neq g_0(c_{v_{\tau},g_0^{-1}})$.

Now, denote
\begin{equation*}
\begin{split}
&U_{\sigma}=\{g\in W\mid g\leq u,\ g^{-1}\leq v_{\sigma}^{-1}\},\\
&U_{\tau}=\{g\in W\mid g\leq u,\ g^{-1}\leq v_{\tau}^{-1}\},
\end{split}
\end{equation*}
then
\begin{equation*}
\begin{split}
&c_{\sigma}=-\dfrac{c_{us_1,g_0}g_0(c_{v_{\sigma},g_0^{-1}})}{\beta}-\sum_{g\in U_{\sigma},\ g\neq g_0}\dfrac{c_{us_1,g}g(c_{v_{\sigma},g^{-1}})}{g\alpha_1},\\
&c_{\tau}=-\dfrac{c_{us_1,g_0}g_0(c_{v_{\tau},g_0^{-1}})}{\beta}-\sum_{g\in U_{\tau},\ g\neq g_0}\dfrac{c_{us_1,g}g(c_{v_{\tau},g^{-1}})}{g\alpha_1}.
\end{split}
\end{equation*}
Suppose
\begin{equation*}
\begin{split}
&-c_{us_1,g_0}=A/B,\ g_0(c_{v_{\sigma},g_0^{-1}})=P_{\sigma}/Q_{\sigma},\ g_0(c_{v_{\tau},g_0^{-1}})=P_{\tau}/Q_{\tau},\\
&-\sum_{g\in U_{\sigma},\ g\neq g_0}\dfrac{c_{us_1,g}g(c_{v_{\sigma},g^{-1}})}{g\alpha_1}=\dfrac{C_{\sigma}}{D_{\sigma}},\\
&-\sum_{g\in U_{\tau},\ g\neq g_0}\dfrac{c_{us_1,g}g(c_{v_{\tau},g^{-1}})}{g\alpha_1}=\dfrac{C_{\tau}}{D_{\tau}}.
\end{split}
\end{equation*}
If $d_{\sigma}=d_{\tau}$, then $c_{\sigma}=c_{\tau}$, too, so
$$\dfrac{A}{B}\cdot\dfrac{P_{\sigma}}{\beta Q_{\sigma}}+\dfrac{C_{\sigma}}{D_{\sigma}}=
\dfrac{A}{B}\cdot\dfrac{P_{\tau}}{\beta Q_{\tau}}+\dfrac{C_{\tau}}{D_{\tau}}.$$
This is equivalent to $$\dfrac{AD_{\sigma}P_{\sigma}+\beta BC_{\sigma}Q_{\sigma}}{\beta BD_{\sigma}Q_{\sigma}}=
\dfrac{AD_{\tau}P_{\tau}+\beta BC_{\tau}Q_{\tau}}{\beta BD_{\tau}Q_{\tau}}.$$

This implies that $$\beta BQ_{\sigma}Q_{\tau}(C_{\sigma}D_{\tau}-C_{\tau}D_{\sigma})=AD_{\sigma}D_{\tau}(P_{\tau}Q_{\sigma}-P_{\sigma}Q_{\tau}).$$ Now, $\beta$ divides neither $A$, nor $D_{\sigma}$, nor $D_{\tau}$, because these non-zero polynomials belong to the subalgebra $S'$ generated by $\alpha_2,\ldots,\alpha_{n-1}$. Since $S=\Cp[\htt]$ is a unique factorization domain, $\beta$ divides $P_{\tau}Q_{\sigma}-P_{\sigma}Q_{\tau}$. But this polynomial belongs to $S'$, thus this polynomial is zero. It follows that $g_0(c_{v_{\sigma},g_0^{-1}})=g_0(c_{v_{\tau},g_0^{-1}})$, a contradiction. Thus, $d_{\sigma}\neq d_{\tau}$. The proof is complete.}

\sst\label{sst:F_4_G_2} In this Subsection, we consider the cases $G_2$ and $F_4$. Actually, in these cases, the proof of Theorem~\ref{mtheo} is by direct computations. For $G_2$, our computations are quite easy. Namely, if $\Phi=G_2$, then there are two fundamental roots $\alpha_1$, $\alpha_2$. The length of the first root is 1, and the length of the second one is $\sqrt{3}$. The angle between $\alpha_1$ and $\alpha_2$ equals $5\pi/6$. Below we list all involutions in the Weyl group of $G_2$ and their Kostant--Kumar polynomials.
\begin{center}
\begin{tabular}{|l@{\hskip\tabcolsep\raise2pt\copy\strutbox\lower2pt\copy\strutbox}|l|}
\hline
$w$&$d_w$\\
\hline\hline
$\id$&$18\alpha_1^5\alpha_2+45\alpha_1^4\alpha_2^2+40\alpha_1^3\alpha_2^3+15\alpha_1^2\alpha_2^4+2\alpha_1\alpha_2^5$\\
\hline
$s_1$&$18\alpha_1^4\alpha_2+45\alpha_1^3\alpha_2^2+40\alpha_1^2\alpha_2^3+15\alpha_1\alpha_2^4+2\alpha_2^5$\\
\hline
$s_1s_2s_1$&$18\alpha_1^3+39\alpha_1^2\alpha_2+27\alpha_1\alpha_2^2+6\alpha_2^3$\\
\hline
$s_1s_2s_1s_2s_1$&$6\alpha_1+4\alpha_2$\\
\hline
$s_2$&$18\alpha_1^5+45\alpha_1^4\alpha_2+40\alpha_1^3\alpha_2^2+15\alpha_1^2\alpha_2^3+2\alpha_1\alpha_2^4$\\
\hline
$s_2s_1s_2$&$18\alpha_1^3+27\alpha_1^2\alpha_2+13\alpha_1\alpha_2^2+2\alpha_2^3$\\
\hline
$s_2s_1s_2s_1s_2$&$4\alpha_1+2\alpha_2$\\
\hline
$s_2s_1s_2s_1s_2s_1$&$1$\\
\hline
\end{tabular}
\end{center}

For $F_4$, our computations are rather complicated. Let $\alpha_1$, $\alpha_2$, $\alpha_3$, $\alpha_4$ be fundamental roots (see \cite{Bourbaki} for the details). For instance, if $w=s_1s_2s_3s_4s_2s_3s_1s_2s_3s_4s_3s_2s_3s_1s_2s_3s_1s_2$, then
\begin{equation*}
\begin{split}
d_w&=4\alpha_1^6+48\alpha_1^5\alpha_2+237\alpha_1^4\alpha_2^2+617\alpha_1^3\alpha_2^3+
894\alpha_1^2\alpha_2^4+684\alpha_1\alpha_2^5+216\alpha_2^6+72\alpha_1^5\alpha_3\\
&+712\alpha_1^4\alpha_2\alpha_3+
2782\alpha_1^3\alpha_2^2\alpha_3+5374\alpha_1^2\alpha_2^3\alpha_3+5136\alpha_1\alpha_2^4\alpha_3+
1944\alpha_2^5\alpha_3+532\alpha_1^4\alpha_3^2\\
&+4160\alpha_1^3\alpha_2\alpha_3^2+12053\alpha_1^2\alpha_2^2\alpha_3^2+
15349\alpha_1\alpha_2^3\alpha_3^2+7254\alpha_2^4\alpha_3^2+2064\alpha_1^3\alpha_3^3+11960\alpha_1^2\alpha_2\alpha_3^3\\
&+22832\alpha_1\alpha_2^2\alpha_3^3+14372\alpha_2^3\alpha_3^3+4432\alpha_1^2\alpha_3^4+
16912\alpha_1\alpha_2\alpha_3^4+15952\alpha_2^2\alpha_3^4+4992\alpha_1\alpha_3^5\\
&+9408\alpha_2\alpha_3^5+2304\alpha_3^6+48\alpha_1^5\alpha_4+476\alpha_1^4\alpha_2\alpha_4
+1862\alpha_1^3\alpha_2^2\alpha_4+3596\alpha_1^2\alpha_2^3\alpha_4\\
&+3432\alpha_1\alpha_2^4\alpha_4+1296\alpha_2^5\alpha_4+712\alpha_1^4\alpha_3\alpha_4
+5568\alpha_1^3\alpha_2\alpha_3\alpha_4+16118\alpha_1^2\alpha_2^2\alpha_3\alpha_4\\
&+20490\alpha_1\alpha_2^3\alpha_3\alpha_4+
9660\alpha_2^3\alpha_3\alpha_4+4144\alpha_1^3\alpha_3^2\alpha_4+
23976\alpha_1^2\alpha_2\alpha_3^2\alpha_4+45676\alpha_1\alpha_2^2\alpha_3^2\alpha_4\\
&+28678\alpha_2^3\alpha_3^2\alpha_4+11840\alpha_1^2\alpha_3^3\alpha_4+
45072\alpha_1\alpha_2\alpha_3^3\alpha_4+42400\alpha_2^2\alpha_3^3\alpha_4+16616\alpha_1\alpha_3^4\alpha_4\\
&+31228\alpha_2\alpha_3^4\alpha_4+9168\alpha_3^5\alpha_4+236\alpha_1^4\alpha_4^2+1844\alpha_1^3\alpha_2\alpha_4^2+
5330\alpha_1^2\alpha_2^2\alpha_4^2+6762\alpha_1\alpha_2^3\alpha_4^2\\
&+3180\alpha_2^4\alpha_4^2+2744\alpha_1^3\alpha_3\alpha_4^2+15884\alpha_1^2\alpha_2\alpha_3\alpha_4^2+30114\alpha_1\alpha_2^2\alpha_3\alpha_4^2+
18858\alpha_2^3\alpha_3\alpha_4^2\\
&+11728\alpha_1^2\alpha_3^2\alpha_4^2+45530\alpha_1\alpha_2\alpha_3^2\alpha_4^2+41776\alpha_2^2\alpha_3^2\alpha_4^2+21868\alpha_1\alpha_3^3\alpha_4^2+
40982\alpha_2\alpha_3^3\alpha_4^2\\
&+15024\alpha_3^4\alpha_4^2+600\alpha_1^3\alpha_4^3+3456\alpha_1^2\alpha_2\alpha_4^3+6552\alpha_1\alpha_2^2\alpha_4^3
+4092\alpha_2^3\alpha_4^3+5112\alpha_1^2\alpha_3\alpha_4^3\\
&+19356\alpha_1\alpha_2\alpha_3\alpha_4^3+18108\alpha_2^2\alpha_3\alpha_4^3
+14244\alpha_1\alpha_3^2\alpha_4^3+26616\alpha_2\alpha_3^2\alpha_4^3+12996\alpha_3^3\alpha_4^3\\
&+828\alpha_1^2\alpha_4^4+3126\alpha_1\alpha_3\alpha_4^4+2916\alpha_2^2\alpha_4^4+
4596\alpha_1\alpha_3\alpha_4^4+8562\alpha_2\alpha_3\alpha_4^4+6264\alpha_3^2\alpha_4^4\\
&+588\alpha_1\alpha_4^5+1092\alpha_2\alpha_4^5+1596\alpha_3\alpha_4^5+168\alpha_4^6.
\end{split}
\end{equation*}
Nevertheless, using the system of computer algebra \texttt{SAGE} \cite{SAGE}, we checked that the Kostant--Kumar polynomials for all 139 involutions in the Weyl group of type $F_4$ are distinct. The listing of our\linebreak program and the complete list of Kostant--Kumar polynomials for involutions are available at\linebreak \texttt{http://algeom.samsu.ru/text/staff-Eliseev.html}. Thus, the proof of Theorem~\ref{mtheo} is complete.

\sect{Concluding remarks}\label{sect:conclusion}

\sst It was conjectured in~\cite{EliseevPanov} that $C_w$ coincides with $C_{w^{-1}}$ for any $w\in W$. The proof of this fact is straightforward, see~\cite{Bochkarev}. On the the other hand, the fact that $d_w=d_{w^{-1}}$ can be easily proved in a purely combinatorial way.

\propp{Let $w \in W$. Let $w^{-1}$ be its inversed. Then $d_{w}= d_{w^{-1}}$.\label{prop:d_w_w_inversed}}{Fix a reduced expression $$w=s_{i_1}s_{i_2}\ldots s_{i_l}.$$
Recall that (see (\ref{formula:d_and_c}) and \cite{KostantKumar2}) $$d_{vw_0,ww_0}=(-1)^{l(w)-l(v)}\cdot c_{w,v}\cdot\prod_{\alpha\in\Phi^+}\alpha.$$ If $v=\id$, then $$d_{w}=d_{w_0,ww_0}=(-1)^{l(w)}\cdot c_{w}\cdot\prod_{\alpha\in\Phi^+}\alpha.$$
Since $l(w)=l(w^{-1})$, $d_{w}= d_{w^{-1}}$ is equivalent to $c_{w}= c_{w^{-1}}$.

By definition, \begin{equation*}
c_w=c_{w,\id}=(-1)^{l(w)}\cdot\sum\dfrac{1}{s_{i_1}^{\epsi_1}\alpha_{i_1}}\cdot\dfrac{1}{s_{i_1}^{\epsi_1}
s_{i_2}^{\epsi_2}\alpha_{i_2}}\cdot\ldots\cdot\dfrac{1}{s_{i_1}^{\epsi_1}\ldots s_{i_l}^{\epsi_l}\alpha_{i_l}},
\end{equation*}
where the sum is taken over all sequences $(\epsi_1,\ldots,\epsi_l)$, $\epsi_i\in\{0,1\}$, such that $s_{i_1}^{\epsi_1}\ldots s_{i_l}^{\epsi_l}=\id$. (Recall that the element $c_{w}$ depends only on $w$, not on the choice of a reduced expression.) We claim that the expressions for $c_{w}$ and $c_{w^{-1}}$ are literally the same (up to order of summands).

Indeed, let $s_{i_1}s_{i_2}\ldots s_{i_l}$ be a reduced expression of $w$, then $s_{i_l}s_{i_{l-1}}\ldots s_{i_1}$ is a reduced expression of $w^{-1}$. Now, consider a sequence $(p_1,\ldots,p_l)$, $p_i\in\{0,1\}$, such that $s_{i_1}^{p_1}\ldots s_{i_l}^{p_l}=\id$. Denote by $k$ the number of units in this sequence. Suppose that $p_{j_1}=p_{j_2}=\ldots=p_{j_k}=1$. Then the summand in the sum for $c_w$ corresponding to $(p_1,\ldots,p_l)$ has the form \begin{equation*}
\begin{split}
&\dfrac{1}{\alpha_{i_1}}\cdot\dfrac{1}{\alpha_{i_2}}\cdot\ldots\cdot\dfrac{1}{s_{i_{j_1}}\alpha_{i_{j_1}}}\cdot\dfrac{1}{s_{i_{j_1}}\alpha_{i_{{j_1}+1}}
}\cdot\ldots\\
&\times\dfrac{1}{s_{i_{j_1}}s_{i_{j_2}}\alpha_{i_{j_2}}}\cdot\dfrac{1}{s_{i_{j_1}}s_{i_{j_2}}\alpha_{i_{{j_2}+1}}
}\cdot\ldots\\
&\times\dfrac{1}{s_{i_{j_1}}s_{i_{j_2}}\ldots s_{i_{j_{k-1}}}\alpha_{i_{j_{k-1}}}}\cdot\dfrac{1}{s_{i_{j_1}}s_{i_{j_2}}\ldots s_{i_{j_{k-1}}}\alpha_{i_{j_{k-1}+1}}
}\cdot\ldots\\
&\times\dfrac{1}{s_{i_{j_1}}s_{i_{j_2}}\ldots s_{i_{j_{k-1}}}\alpha_{i_{j_k-1}}}\cdot\dfrac{1}{\alpha_{i_{j_k}}}\cdot\dfrac{1}{\alpha_{i_{j_{k+1}}}}\cdot\ldots\cdot\dfrac{1}{\alpha_{i_{l}}}.
\end{split}
\end{equation*}
Consider the summand in the sum for $c_{w^{-1}}$ corresponding to the sequence $(p_l,\ldots,p_1)$. (It is clear that $s_{i_l}^{p_l}\ldots s_{i_1}^{p_1}=\id$, because $s_{i_1}^{p_1}\ldots s_{i_l}^{p_l}=\id$.) In this sequence, units are situated on the places $l-j_1+1$, $l-j_2+1$, $\ldots$, $l-j_k+1$. Denote $s_{i_j}'=s_{i_{l-i+1}}$ and $\alpha_{i_j}'=\alpha_{i_{l-i+1}}$.

Let $1\leq t\leq k$. Consider the factor of the denominator of this summand in $c_{w^{-1}}$ of the form $$\dfrac{1}{s_{i_{l-j_k+1}}'s_{i_{l-j_{k-1}+1}}'\ldots s_{i_{l-j_t+1}}'\alpha_{i_{l-j_t+2}}'}\cdot\ldots\cdot\dfrac{1}{s_{i_{l-j_k+1}}'s_{i_{l-j_{k-1}+1}}'\ldots s_{i_{l-j_t+1}}'\alpha_{i_{l-j_{t-1}}}'}.$$ Since $$s_{i_{l-j_k+1}}'s_{i_{l-j_{k-1}+1}}'\ldots s_{i_{l-j_t+1}}'s_{i_{l-j_{t-1}+1}}'\ldots s_{i_{l-j_1+1}}'=\id,$$ we obtain $$s_{i_{l-j_k+1}}'s_{i_{l-j_{k-1}+1}}'\ldots s_{i_{l-j_t+1}}'=s_{i_{l-j_1+1}}'\ldots s_{i_{l-j_{t-1}+1}}'=s_{j_1}\ldots s_{j_{t-1}}.$$ But the denominator of the summand for $c_w$ corresponding to $(p_1,\ldots,p_l)$ has the factor $$\dfrac{1}{s_{i_{j_1}}\ldots s_{i_{j_{t-1}}}\alpha_{i_{j_{t-1}+1}}}\cdot\ldots\dfrac{1}{s_{i_{j_1}}\ldots s_{i_{j_{t-1}}}\alpha_{i_{j_t-1}}}.$$
Since $$\dfrac{1}{s_{i_{l-j_k+1}}'s_{i_{l-j_{k-1}+1}}'\ldots s_{i_{l-j_t+1}}'\alpha_{i_{l-j_t+1}}'}=
\dfrac{1}{s_{i_{j_1}}\ldots s_{i_{j_{t-1}}}\alpha_{i_{j_t}}}=-\dfrac{1}{s_{i_{j_1}}\ldots s_{i_{j_{t-1}}}s_{i_{j_t}}\alpha_{i_{j_t}}}$$ and $k$ is even, we conclude that the summand in $c_w$ corresponding to $(p_1,\ldots,p_l)$ is equal to the summand in $c_{w^{-1}}$ corresponding to $(p_l,\ldots,p_1)$. The result follows.}

\sst\label{sst:coadjoint} In this Subsection, we briefly describe interaction of geometry of tangent cones with coadjoint orbits for the case $G=\GL_n(\Cp)$ or $\SL_n(\Cp)$ (according to A.A. Kirillov's orbit method \cite{Kirillov1}, \cite{Kirillov2}, coadjoint orbits play the key role in representation theory of unipotent radical $U$ of the group $B$). Here $U$ is the unitriangular group, i.e., the group of all upper-triangular matrices with $1$'s on the diagonal; its Lie algebra $\nt$ consists of all upper-triangular matrices with zeroes on the diagonal. The groups $B$ and $U$ act on $\nt$ by the adjoint action (i.e., by conjugation). The dual action on the dual space~$\nt^*$ is called \emph{coadjoint}. Recall that we denote by $\gt$, $\bt$ the Lie algebras of the groups $G$, $B$ respectively.

The tangent space $T_p\Fo$ to the flag variety $\Fo=G/B$ can be naturally identified with $\gt/\bt$. Using the Killing form on $\gt$, one can identify the latter space with $\nt^*$. Thus, the tangent cones $C_w$, $w\in S_n$, are subschemes of $T_pX_w\subseteq T_p\Fo\cong\nt^*$. Further, the action of $B$ on $\Fo$ by conjugation induces the action of $B$ on the tangent space $T_p\Fo$. In fact, this action coincides with the coadjoint action of $B$ on $\nt^*$ under the identification $T_p\Fo\cong\nt^*$. Each tangent cone is $B$-invariant, i.e., is a union of coadjoint orbits.

On the other hand, to each involution $w\in S_n$ one can assign the coadjoint orbit $\Omega_w\subseteq\nt^*$ of $B$ by the following rule. Consider the standard basis of $\nt$ consisting of matrix units. Denote by $f_w$ the element of $\nt^*$ equal to the sum of covectors $e_{j,i}^*$, $j < i$, such that $w(i)=j$. It is easy to see that $\Omega_w\subseteq C_w$, so $\overline{\Omega}_w\subseteq C_w$. Computations in \cite{EliseevPanov} and \cite{Ignatyev2} show that if $n\leq5$, then the closure $\overline{\Omega}_w$ coincides with the tangent cone $C_w$ for all involutions $w\in S_n$. Hence we have the following conjecture \cite[Conjecture 1.11]{Ignatyev2}: $\overline{\Omega}_w=C_w$ for all involutions $w\in S_n$. See \cite[Section 4]{Ignatyev2} for the dimension of~$\Omega_w$, conjectural description of $\overline{\Omega}_w$ and further details.

\end{document}